\documentclass[12pt]{article}
\usepackage{amssymb,latexsym,amsmath,amsfonts}

\newtheorem{theorem}{Theorem}[section]

\newtheorem{corollary}[theorem]{Corollary}

\newtheorem{conjecture}{Conjecture}[section]

\newtheorem{problem}[theorem]{Problem}

\newcounter{parag}[section]
\newcounter{remark}[theorem]

\newenvironment{parag}[1]{
\noindent{\bf\boldmath\stepcounter{parag}\arabic{section}.\arabic{parag} #1.}}

\newenvironment{proof}{\noindent{\it Proof:}}{$\square$}

\newenvironment{remark}{\noindent {\it Remark} \stepcounter{remark}\arabic{section}.\arabic{theorem}.\arabic{remark}}{\vskip0.4cm}
\newenvironment{conjecturemark}{\noindent {\it Remark} \stepcounter{remark}\arabic{section}.\arabic{conjecture}.\arabic{remark}}{\vskip0.4cm}

\def\beq{\begin{eqnarray}}
\def\eeq{\end{eqnarray}}
\def\bes{\begin{eqnarray*}}
\def\ees{\end{eqnarray*}}

\def\C{\mathbb{C}}
\def\M{{\mathcal{M}}^d}

\def\F{\mathbb{F}}

\def\Z{\mathbb{Z}}
\def\N{\mathcal{N}^d}

\newcommand{\nc}{\newcommand}
\nc{\op}[1]{\mathop{\mathchoice{\mbox{\rm #1}}{\mbox{\rm #1}}
{\mbox{\rm \scriptsize #1}}{\mbox{\rm \tiny #1}}}\nolimits}
\nc{\al}{\alpha}

\nc{\ep}{\varepsilon}
\nc{\ga}{\gamma}
\nc{\Ga}{\Gamma}
\nc{\la}{\lambda}
\nc{\La}{\Lambda}
\nc{\si}{\sigma}
\nc{\Sig}{{\Gamma}}
\nc{\Om}{\Omega}
\nc{\om}{\omega}
\nc{\SL}[1]{{{\rm SL(}#1{\rm )}}}
\nc{\GL}[1]{{{\rm GL(}#1{\rm )}}}
\nc{\PGL}[1]{{{\rm PGL(}#1{\rm )}}}
\nc{\U}[1]{{{\rm U(}#1{\rm )}}}
\nc{\cpt}{{\op{cpt}}}
\nc{\Dol}{{\op{Dol}}}
\nc{\DR}{{\op{DR}}}
\nc{\B}{{\op{B}}}
\nc{\Triv}{\op{Triv}}
\nc{\Hod}{{\op{Hod}}}
\nc{\Est}{E_{\op{st}}}
\nc{\Hst}{H_{\op{st}}}
\nc{\Left}[1]{\hbox{$\left#1\vbox to
    10.5pt{}\right.\nulldelimiterspace=0pt \mathsurround=0pt$}}
\nc{\Right}[1]{\hbox{$\left.\vbox to
    10.5pt{}\right#1\nulldelimiterspace=0pt \mathsurround=0pt$}}
\nc{\LEFT}[1]{\hbox{$\left#1\vbox to
    15.5pt{}\right.\nulldelimiterspace=0pt \mathsurround=0pt$}}
\nc{\RIGHT}[1]{\hbox{$\left.\vbox to
    15.5pt{}\right#1\nulldelimiterspace=0pt \mathsurround=0pt$}}

\nc{\bee}{{\bf E}}
\nc{\bphi}{{\bf \Phi}}

\begin{document}
\thispagestyle{empty}

\title{Mirror symmetry and Langlands duality in the non-Abelian Hodge theory of a curve}

\author{
 Tam\'as Hausel
\\ {\it Department of Mathematics}
\\ {\it University of Texas at Austin}
\\ {\it Austin TX 78712, USA}
\\{\tt hausel@math.utexas.edu}}
\maketitle

\begin{abstract} 
This is a survey of results and 
conjectures 
on mirror symmetry phenomena in the non-Abelian Hodge theory of a 
curve. We start with the 
conjecture of Hausel--Thaddeus which claims that certain 
Hodge numbers of moduli spaces of flat $SL(n,\C)$ and $PGL(n,\C)$-connections
on a smooth projective algebraic curve agree. 
We then change our point
of view in the non-Abelian 
Hodge theory of the curve, and concentrate on the $SL(n,\C)$ and 
$PGL(n,\C)$ character varieties of 
the curve. Here we discuss a recent
conjecture of Hausel--Rodriguez-Villegas which claims, analogously to
the above conjecture, that 
certain Hodge numbers of these character varieties also agree. We 
explain that for Hodge numbers of character varieties one can use
arithmetic methods, and thus we 
end up explicitly calculating, in terms of Verlinde-type formulas,
the number of representations
of the fundamental group into the finite groups $SL(n,\F_q)$ and
$PGL(n,\F_q)$, by using the character tables of these finite groups of Lie type.  
Finally we explain a conjecture  
which enhances the previous result, and  
gives a simple formula for the  mixed Hodge polynomials, and in
particular for the Poincar\'e 
polynomials of these character varieties, and detail the relationship
to results of Hitchin, Gothen, Garsia--Haiman and Earl--Kirwan. One
consequence of this conjecture is a curious Poincar\'e 
duality type of symmetry, which leads to a conjecture, similar to
Faber's conjecture on the moduli space of curves,  
about a strong Hard Lefschetz theorem for the character variety, which
can be considered as a generalization 
of both the Alvis--Curtis duality in
the representation theory of finite groups of Lie type and a 
recent result of the author on the quaternionic geometry of matroids.

\end{abstract}

\section{Introduction} Non-Abelian Hodge theory (\cite{hitchin}, \cite{simpson}) 
of a genus $g$ smooth complex projective curve $C$
studies three moduli spaces  attached to $C$ and a 
reductive complex algebraic group $G$, which in this paper
will be either $GL(n,\C)$ or $SL(n,\C)$ or $PGL(n,\C)$. They are 
$\M_{Dol}(G)$, the moduli space of semistable 
$G$-Higgs bundles on $C$; $\M_{DR}(G)$,
the moduli space of flat $G$-connections on $C$ and $\M_B(G)$ the character variety, i.e. the 
moduli space of representations of $\pi_1(C)$ into $G$ modulo conjugation. 
With some assumptions these  moduli spaces are smooth varieties (or orbifolds
when $G=PGL(n,\C)$) with the 
underlying differentiable manifolds canonically identified, which carries a 
natural hyperk\"ahler metric. 

The cohomology of this underlying manifold
has mostly been studied from the perspective of $\M_{Dol}(G)$. Using a natural circle action
on it \cite{hitchin} and \cite{gothen} calculated the Poincar\'e polynomials for $G=SL(2,\C)$ and 
$G=SL(3,\C)$ respectively; while \cite{HT1} and \cite{markman} found a simple set of generators for the cohomology ring for $G=PGL(2,\C)$ and $G=PGL(n,\C)$ respectively. The paper $\cite{HT2}$ then 
calculated the cohomology ring explicitly for $G=PGL(2,\C)$. The techniques used in these papers
do not seem to generalize easily to higher $n$. 

A new perspective for these investigations on the cohomology of $\M_{Dol}(G)$ and $\M_{DR}(G)$ 
was introduced in \cite{HT3} and \cite{HT4}. It was shown there that the hyperk\"ahler 
metrics and the Hitchin systems \cite{hitchin2} 
for $\M_{DR}(G)$ and $\M_{DR}(G^L)$, with $G=SL(n,\C)$ and Langlands dual $G^L=PGL(n,\C)$ provide 
the geometrical setup suggested in \cite{SYZ} as a criteria for mirror symmetry.  
Based on this observation \cite{HT4} conjectured that a version of the topological mirror symmetry
also holds, i.e. that certain Hodge numbers for $\M_{DR}(G)$ and
$\M_{DR}(G^L)$ agree. Using the above mentioned results of \cite{hitchin} and \cite{gothen} this conjecture
was checked for $G=SL(2,\C), SL(3,\C)$. This mirror symmetry conjecture motivates the study of
not just the cohomology but the mixed Hodge structure on the cohomology of the spaces $\M_{DR}(G)$,
$\M_{Dol}(G)$ and $\M_{B}(G)$. While it was shown in \cite{HT4} that the mixed Hodge structure
of $\M_{Dol}(G)$ and $\M_{DR}(G)$ agree, and can be shown to be pure as in Theorem~\ref{pure}
or \cite{mehta}, the mixed Hodge structure on $\M_{B}(G)$ has not been
studied until very recently.  

An important theme of this survey paper is in fact the mixed Hodge structure on 
the  character variety 
$\M_{B}(G)$ or alternatively   the three variable 
polynomial $H(x,y,t)$ the so-called mixed Hodge polynomial or shortly 
$H$-polynomial which encodes the dimensions of the graded 
pieces of the mixed Hodge 
structure on
$\M_{B}(G)$.
 In a recent project \cite{HRV} an arithmetic method was used 
to calculate the  $E$-polynomial (where the E-polynomial $E(x,y)$ for a smooth variety is defined as $x^ny^nH(1/x,1/y,-1)$, where $n$ is the 
complex dimension of the variety) 
of $\M_{B}(G)$. The idea of \cite{HRV} is to count
the rational points of $\M_B(G(\F_q))$, the variety $\M_{B}(G)$ over the finite field $\F_q$, 
where $q$ is a prime power. This count then is possible due to a result of 
\cite{mednykh}, and the resulting formula, which resembles the famous Verlinde formula \cite{verlinde}, 
is a simple sum over irreducible representations of the finite group of Lie type $G(\F_q)$.
Thus the representation theory behind the $E$-polynomial of the character variety is that of the 
finite groups of Lie type, which could be considered as an 
analogue
of Nakajima's principle \cite{nakajima}, which states 
that the representation theory of Kac--Moody algebras are 
encoded in the cohomology of the (hyperk\"ahler) quiver varieties.  

The shape of the $E$-polynomials of the various character varieties 
then made us conjecture \cite{HRV} that the mirror symmetry conjecture 
also holds for the pair $\M_{B}(G)$ and $\M_{B}(G^L)$ in the case of
$G=SL(n,\C)$ at 
least. 
Due to our
improved ability to calculate these Hodge numbers via this number
theoretical method, 
we could 
check this conjecture in the cases when $n$ is $4$ or a prime. As the
two mirror symmetry 
conjectures of \cite{HT4} and \cite{HRV} are equivalent on the level
of Euler characteristic, we 
get a proof of the original mirror symmetry conjecture of \cite{HT4} 
on the level of Euler characteristic in these cases. 

Perhaps even more interestingly \cite{HRV} achieves explicit formulas, in terms of 
a simple 
generating function, for the $E$-polynomials
of the character variety $\M_B(GL(n,\C))$. In particular it can be deduced 
from this that the Euler 
characteristic  of $\M_B(PGL(n,\C))$ is $\mu(n) n^{2g-3}$, where $\mu$ is 
the fundamental number theoretic function: the M\"obius function, i.e. the sum of all 
primitive $n$th roots of unity. This result, which could not be 
obtained with other methods, in itself hints 
at an interesting link between number theory and the 
topology of our hyperk\"ahler manifolds. Another consequence of our
formula, is that these $E$-polynomials turn out to be palindromic,
i.e. satisfy an unexpected Poincar\'e duality-type of symmetry. In
fact this symmetry can be traced back to the Alvis--Curtis duality
\cite{alvis,curtis} in the representation theory of finite
groups of Lie type.

Then we present a deformation of the formula for the $E$-polynomial of the character variety 
$\M_B(PGL(n,\C))$, which conjecturally \cite{HRV} should agree with the 
$H$-polynomial.   
Moreover 
we will later 
modify this formula to obtain what conjecturally should be the $H$-polynomial of 
the Higgs moduli space $P\M_{Dol}(GL(n,\C))$. We also explain how,
using, 
as a guide, our mirror symmetry conjectures, 
one could get conjectures for the corresponding 
$H$-polynomials for the varieties associated to
$SL(n,\C)$. 

These conjectures imply a conjecture on the Poincar\'e 
polynomials (where the Poincar\'e polynomial is obtained from the
$H$-polynomial as $H(1,1,t)$) of our manifolds $\M_{Dol}(PGL(n,\C))$. 
This conjecture is similar in flavour to
Lusztig's conjecture \cite{lusztig} on the Poincar\'e polynomials of  Nakajima's
quiver varieties, which is also a hyperk\"ahler manifold, similar to
the Higgs moduli space $\M_{Dol}(G)$. We should also mention Zagier's
\cite{zagier} formula for the Poincar\'e polynomial of the moduli space $\N$ of stable
bundles (the ``K\"ahler version'' of $\M_{Dol}(SL(n,\C))$), 
where the formula is
a similar sum, but is parametrized by ordered partitions of $n$. 

We discuss in detail many checks on these conjectures, by showing how our 
conjectures imply results obtained by Hitchin \cite{hitchin}, Gothen
\cite{gothen} and 
Earl--Kirwan 
\cite{earl-kirwan}. Already the combinatorics of these formulas are
non-trivial, and  surprisingly 
the calculus of Garsia--Haiman \cite{garsia-haiman} needs to be used to
check the  conjecture when $g=0$. 

Since a curious Poincar\'e  duality type of symmetry is satisfied for
the conjectured Hodge numbers of  $\M_{B}(PGL(n,\C))$, 
we also discuss the conjecture that a certain version of the Hard
Lefschetz theorem is satisfied 
for our non-compact varieties. This is then explained to be a
generalization  of a result in 
\cite{hausel3} on the quaternionic geometry of matroids, and an analogue of the Faber conjecture 
\cite{faber}
on the moduli space of curves.  

\paragraph{\bf Acknowledgment} This survey paper is based on the author's talk at the 
``Geometric Methods in Algebra and Number Theory'' conference at the 
University of Miami in December 2003. I would 
like to thank the organizers for the invitation 
and for the memorable conference. 
Also most of the results and conjectures 
surveyed here have recently been obtained in joint projects with 
Michael Thaddeus \cite{HT3,HT4} and 
with Fernando Rodriguez-Villegas \cite{HRV}. 
The research described in this paper was partially supported by the 
NSF grant DMS-0305505.

\section{Abelian and non-Abelian Hodge theory}

This section gives some basic definitions on Abelian and 
non-Abelian Hodge theory which will be used in the paper later. 
For details on them consult the sources indicated below. 

\vskip.5cm

\begin{parag}{Hodge--De Rham theory} Fix a 
smooth complex algebraic variety $M$. There are various 
cohomology theories
which associate a graded anti-commutative ring to the variety $M$. First we 
consider the singular, or Betti, cohomology $H^*_B(M,\C)$ of $M$ 
with complex coefficients. 
This in fact can be defined for any reasonable 
topological space. 
The dimension of $H^k_B(M,\C)$ is called the $k$-th Betti number 
and denoted $b_k(M)$. The Poincar\'e polynomial 
is then formed from these numbers as coefficients: $$P(t;M)=\sum_k b_k(M) t^k.$$ 

Next we consider the De Rham cohomology 
$H^*_{DR}(M,\C)$, which is the space of closed differential forms modulo 
exact ones. This can be defined on any differentiable manifold. The De Rham 
theorem then shows that these two cohomologies are naturally isomorphic:
\beq \label{derham} H^*_B(M,\C)\cong H^*_{DR}(M,\C).\eeq Now we assume that our variety is 
projective. Then we have the Dolbault cohomology $H^*_{Dol}(M,\C)$, which
is defined as $$H^k_{Dol}(M,\C)=\bigoplus_{p+q=k} H^q(M,\Omega^p_M).$$
The Hodge theorem then implies that there is a natural isomorphism
\beq \label{hodge}H^k_{DR}(M,\C)\cong H^k_{Dol}(M,\C).\eeq The above two isomorphisms
then imply the Hodge decomposition theorem: 
\beq H^k_{B}(M,\C)\cong \bigoplus_{p+q=k} H^{p,q}(M), \label{decomp}\eeq where
$H^{p,q}(M)$  denotes $H^p(M,\Omega^q_M)$. The dimension of $H^{p,q}(M)$ is denoted $h^{p,q}(M)$ and called the Hodge numbers of the variety $M$. From these numbers we form a two variable polynomial the 
Hodge polynomial: $$H(x,y;M):=\sum_{p,q} h^{p,q}(M) x^p y^q.$$

For more details on these cohomology theories see \cite{griffiths-harris}.

\end{parag}
\vskip.5cm
\begin{parag}{Mixed Hodge structures}
Deligne \cite{deligne2} generalized the Hodge decomposition theorem (\ref{decomp}) to any complex variety $M$, not necessarily smooth or projective, by introducing a so called {\em mixed Hodge structure} on $H_B^*(M,\C)$. 
This implies a 
decomposition\footnote{In fact what one gets from a mixed Hodge structure are two filtrations on the cohomology, 
and the decomposition in question
is the associated graded.} $$H_B^k(M,\C)\cong \bigoplus_{p,q} H^{p,q;k}(M), $$ 
where $p+q$ is called the weight of $H^{p,q;k}(M)$.
In the case of a smooth projective variety we have 
$H^{p,q;p+q}(M)=H^{p,q}(M)$, i.e. that the weight of $H^{p,q;k}(M)$ 
is always $k$, this weight is called {\em pure weight}. 
However in general
we could have other weights appear in the mixed Hodge structure of a
complex algebraic variety; indeed we will see such examples later.  
The dimensions of $H^{p,q;k}(M)$ are 
denoted by $h^{p,q;k}(M)$ and are called mixed Hodge numbers. 
From them we form the three variable polynomial 
\beq H(x,y,t;M):=\sum_{p,q,k} h^{p,q;k}(M) x^p y^q t^k.\label{mixedhodge}\eeq

Similarly, Deligne \cite{deligne2} constructs a mixed Hodge structure
on the compactly supported $H^*_{B,cpt}(M,\C)$ singular cohomology
of our complex algebraic variety $M$. This yields the decomposition 
 $$H_{B,cpt}^k(M,\C)\cong \bigoplus_{p,q} H_{cpt}^{p,q;k}(M), $$ and leads
to the compactly supported mixed Hodge numbers $h^{p,q;k}_{cpt}(M)$,
which is defined as the dimension of $H_{cpt}^{p,q;k}(M)$. Then one
can introduce the $e$-numbers $e^{p,q}(M) = \sum_{k} (-1)^k h_{cpt}^{p,q;k}(M)$ from which we get the $E$-polynomial: \begin{eqnarray}E(x,y;M):=\sum_{p,q} e^{p,q}(M) x^p y^q.\label{epolynomial}\end{eqnarray} 
Clearly for a smooth projective variety $E(x,y)=H(-x,-y)$. Moreover for a smooth variety Poincar\'e duality implies 
that $$E(x,y)=(xy)^n H(1/x,1/y,-1),$$ where $n$ is the complex dimension of $M$. 
The significance of the $E$-polynomial is that it is additive 
for decompositions and multiplicative for Zariski locally trivial 
fibrations.  

For more details see the original \cite{deligne2} or \cite{batyrev-dais} 
for more on the $E$-polynomials. 
\end{parag}
\vskip.5cm
\begin{parag}{Stringy cohomology} \label{stringy} Suppose a finite group $\Gamma$ acts on our variety $M$. Then by the naturality of the mixed Hodge structure $\Gamma$ will act on 
$H^{p,q,k}(M)$ and we have $$H^{p,q;k}(M/\Gamma)\cong \left(H^{p,q;k}(M)\right)^\Gamma.$$ However for a Calabi--Yau $M$ and $\Gamma$ preserving the Calabi--Yau structure
 string theorists \cite{vafa, zaslow} introduced different Hodge 
numbers on the Calabi--Yau orbifold $M/\Gamma$: the so-called stringy Hodge numbers, which are the right
Hodge numbers for mirror symmetry. Their mathematical
significance is highlighted by a theorem of Kontsevich \cite{kontsevich} 
that the stringy Hodge numbers agree with the ordinary 
Hodge numbers of any 
crepant resolution.  
Following \cite{batyrev-dais} we can define the stringy
$E$-polynomials:
$$\Est(x,y;M/\Ga) = \sum_{[\gamma]} E(x,y;M^\gamma)^{C(\gamma)}
(xy)^{F(\gamma)}.$$
Here the sum runs over the conjugacy classes of
$\Ga$; $C(\ga)$ is the centralizer of $\ga$; $M^\ga$ is the subvariety
fixed by $\ga$; and $F(\ga)$ is an integer called the {\em fermionic
shift}, which is defined as follows.
The group element $\ga$ has finite order, so it acts on $TM|_{M^\ga}$
as a linear automorphism with eigenvalues $e^{2 \pi i w_1}, \dots,
e^{2 \pi i w_n}$, where each $w_j \in [0,1)$.  Let $F(\ga) = \sum
w_j$; this is an integer since, by hypothesis, $\ga$ acts trivially on
the canonical bundle. 

The last cohomology theory we will need is the stringy cohomology 
of a Calabi--Yau orbifold twisted by a $B$-field. 
Following \cite{hitchin3} we let 
$B\in H^2_\Gamma (M, U(1))$ i.e. an isomorphism class of a $\Ga$-equivariant flat unitary 
gerbe. For any element $\ga\in \Gamma$ this $B$-field 
induces a $C(\ga)$-equivariant local system $[L_{B,\gamma}]\in H^1_{C(\gamma)} (M^\gamma,U(1))$ on the fixed point set $M^\gamma$. Using
this we can twist the stringy $E$-polynomial of $M/\Gamma$ to get:
\beq \Est^B (x,y;M / \Ga) = \sum_{[\ga]}
E(x,y;M^\ga;L_{B,\ga})^{C(\gamma)}(xy)^{F(\ga)}.\label{bfieldstringy}\eeq

For more information on the mathematics of stringy cohomology see \cite{batyrev-dais}, for twisting with a $B$-field see \cite{HT4}. 
\end{parag}
\vskip.5cm
\begin{parag}{Non-Abelian Hodge theory} The starting
point of non-Abelian Hodge theory is the identification of the space $H_B^1(M,\C^\times)$ with the space of homomorphisms from 
$\pi_1(M)\to \C^\times$; the space $H_{DR}^1(M,\C^\times)$ with
algebraic local systems on $M$ and the space $$H_{Dol}(M,\C^\times)\cong
H^1(M,{\cal O}^\times)\oplus H^0(M,\Omega^1)$$ with pairs of 
a holomorphic line bundle and a holomorphic one form.

This then can be generalized to any 
non-Abelian complex reductive group $G$. We define $H_B^1(M,G)$
to be conjugacy classes of representations of $\pi_1(M)\to G$. 
I.e. $$H_B^1(M,G):={\rm Hom}(\pi_1(M),G)// G,$$ which is the affine GIT quotient of the affine variety ${\rm Hom}(\pi_1(M),G)$ by the conjugation action by $G$. This is sometimes called the {\em character variety}. 
The space $H_{DR}^1(M,G)$ can be identified as the moduli space of algebraic $G$-local systems on $M$. Finally $H_{Dol}^1(M,G)$ is defined as the moduli space of certain semistable $G$-Higgs bundles on $M$. We will 
give precise definition in the case of a curve below. 
The identification between  $H_B^1(M,G)$ and $H_{DR}^1(M,G)$, 
which is analogous to the De Rham map (\ref{derham}),  is given by
the Riemann--Hilbert correspondence \cite{deligne1,simpson4}, 
while the identification 
between  $H_{DR}^1(M,G)$ and $H_{Dol}^1(M,G)$, which is analogous to 
the Hodge decomposition (\ref{hodge}), is given in \cite{corlette, simpson2} 
by the theory
of harmonic bundles, which is the non-Abelian generalization of 
Hodge theory. 

For an introduction to non-Abelian Hodge theory see \cite{simpson}, and
\cite{katzarkov-pantev}[Section 3], for more details on the construction of the spaces
appearing in non-Abelian Hodge theory and the maps between them 
see \cite{simpson2,simpson3,simpson4}.
\end{parag}
\vskip.5cm
\begin{parag}{The case of a curve} We now fix a smooth 
projective complex curve $C$ of genus $g$ 
and specify our spaces in the case when $M= C$ and $G=GL(n,\C)$. 
According to the above definitions we have:
\bes{\cal M}_B(GL(n,\C))&:=&H^1_B(C,GL(n,\C))= \\ &=& 
\{A_1,B_1,\dots,A_g,B_g\in GL(n,\C) | [A_1,B_1]\cdot \dots \cdot 
[A_g,B_g]=Id \}// GL(n,\C).\ees
There is a natural way to twist these varieties, and because they 
will be needed for $PGL(n,\C)$ we introduce these twisted varieties 
here. For an integer $d$ we also consider:
\bes\M_B(GL(n,\C)):= \{A_1,B_1,\dots,A_g,B_g\in GL(n,\C) | [A_1,B_1]
\cdot \dots \cdot [A_g,B_g]=e^{\frac{2\pi i d}{n}}Id \}// GL(n,\C).
\ees

The De Rham space looks like
\bes{\mathcal M}_{DR}(GL(n,\C))&:=&H^1_{DR}(C,GL(n,\C))= \\ &=& 
\{\mbox{moduli
space of flat $GL(n,\C)$-connections on $C$}\}\ees
and in the twisted case we need to fix a point $p\in C$, and define
\bes {\mathcal M}^d_{DR}(GL(n,\C)) := \left\{ \begin{array}{c} \mbox{ 
moduli space of flat $GL(n,\C)$-connections on $C\setminus\{ p \}$},
\\
\mbox{with holonomy $e^{\frac{2\pi i d}{n}}Id$ around $p$}
\end{array}\right\}.\ees

Finally the Dolbeault spaces are:
\bes{\mathcal M}_{Dol}(GL(n,\C))&:=&H^1_{Dol}(C,GL(n,\C))= \\ &=& 
\{\mbox{moduli
space of semistable  rank $n$ degree $0$ Higgs bundles on $C$}\},\ees
where a rank $n$ Higgs bundle is a pair $(E,\phi)$ of a rank $n$ 
algebraic vector bundle $E$ on $C$, with degree $0$ and Higgs field 
$\phi\in H^0(C,K{\rm End} E)$. A Higgs bundle is called 
semistable
if for any Higgs subbundle $(F,\psi)$ (i.e. a subbundle with compatible Higgs fields) we have $\frac{{\rm deg}(F)}{{\rm rank}(F)}\leq \frac{{\rm deg}(E)}{{\rm rank}(E)}=0$.
The twisted version of ${\mathcal M}_{Dol}(GL(n,\C))$ is defined:
\bes{\mathcal M}^d_{Dol}(GL(n,\C)):=\{\mbox{ moduli
space of semistable  rank $n$ degree $d$ Higgs bundles on $C$}\}.\ees

The varieties defined above for $GL(n,\C)$ are all of dimension $n^2(2g-2)+2$.
The Betti space is affine, while the De Rham space is
analytically isomorphic, via the Riemann--Hilbert correspondence, 
to the Betti space but not algebraically, so the De Rham space is
a Stein manifold as a complex manifold 
but not an affine variety as an algebraic variety. Finally the Dolbeault 
space is a quasi projective variety with large projective subvarieties. 

From now on we will consider only the case when $(n,d)=1$; so we 
fix such a $d$. In this 
case the corresponding twisted spaces are additionally smooth, have a diffeomorphic underlying manifold $\M(GL(n,\C))$ which carries 
a complete hyperk\"ahler metric \cite{hitchin}. 
The complex 
structures of $\M_{Dol}(GL(n,\C))$ and $\M_{DR}(GL(n,\C))$ 
appear in the hyperk\"ahler structure.  

We started this subsection by determining these spaces for $GL(1,\C)
\cong \C^\times$.
By the identifications explained there we see that \beq {\M}_B(GL(1,\C))
&\cong& (\C^\times)^{2g},\nonumber \\  
\M_{Dol}(GL(1,\C))&\cong& T^* Jac^d(C)\label{gl1}\eeq and
$\M_{DR}$ is a certain affine bundle over $Jac^d(C)$. Interestingly for $d=0$
they are all 
algebraic groups and they act on the corresponding
spaces for $GL(n,\C)$ and any $d$ by tensorization. 

We can consider the 
map  $$\begin{array}{cccc} \lambda_{Dol}:&\M_{Dol}(GL(n,\C))&\to& \M_{Dol}(GL(1,\C))\\ &(E,\Phi)&\mapsto &({\det}(E),{\rm tr}(\phi)).\end{array}$$ 
The fibres of this map can be shown to be isomorphic using the 
above tensorization action. So up to isomorphism it is irrelevant which fibre we take, but we usually take a point $(\Lambda,0)\in  \M_{Dol}(GL(1,\C))$ and define $$\M_{Dol}(SL(n,\C)):=\lambda_{Dol}^{-1}((\Lambda,0)),$$ for the other two spaces we have:

\bes {\mathcal M}^d_{DR}(SL(n,\C)) = \left\{ \begin{array}{c} \mbox{ 
moduli space of flat $SL(n,\C)$-connections on $C\setminus\{ p \}$}
\\
\mbox{with holonomy $e^{\frac{2\pi i d}{n}}Id$ around $p$}
\end{array}\right\},\ees
and
\bes\M_B(SL(n,\C))= \{A_1,B_1,\dots,A_g,B_g\in SL(n,\C) | [A_1,B_1]
\cdot \dots \cdot [A_g,B_g]=e^{\frac{2\pi i d}{n}}Id \}// SL(n,\C).
\ees

 The varieties $\M_{B}(SL(n,\C))$,   $\M_{DR}(SL(n,\C))$ and $\M_{DR}(SL(n,\C))$ are smooth of dimension $(n^2-1)(2g-2)$, with diffeomorphic 
underlying manifold $\M(SL(n,\C))$. The Betti space
is affine, and the Betti and De Rham spaces are again analytically, but not
algebraically, isomorphic.  

Moreover we see that a finite subgroup namely 
$Jac[n]\cong \Z_n^{2g}\subset {\mathcal M}_{Dol}(GL(1,\C))$ preserves the 
fibration $\lambda_{Dol}$ and thus acts on $\M_{Dol}(SL(n,\C))$. The quotient then is: 
\bes \M_{Dol}(PGL(n,\C)):=\M_{Dol}(SL(n,\C))/{\rm Jac}[n]\ees
and similarly
\bes \M_{DR}(PGL(n,\C)):=\M_{DR}(SL(n,\C))/{\rm Jac}[n], \ees
and \bes \M_{B}(PGL(n,\C))=\M_{B}(SL(n,\C))/\Z_n^{2g}.\ees This
shows that the spaces  $\M_{B}(PGL(n,\C))$, $\M_{DR}(PGL(n,\C))$ and 
$\M_{Dol}(PGL(n,\C))$ are hyperk\"ahler orbifolds of dimension 
$(n^2-1)(2g-2)$. As they are orbifolds 
we can talk about their stringy mixed Hodge numbers as defined above
in 2.\ref{stringy}. Moreover there are natural orbifold $B$-fields on 
them, which we construct now. First we consider 
a  universal Higgs pair
$(\bee, \bphi)$ on $\M_\Dol(SL(n,\C))\times C$, it exists
because $(d,n)=1$. 
Restrict $\bee$ to $\M_\Dol\times \{p\}$ to get the vector 
bundle $\bee_p$ on $\M_{Dol}(SL(n,\C))$.  
Now we can consider the projective bundle ${\mathbb P}\bee_p$ of 
$\bee_p$ which is a $PGL(n,\C)$-bundle. The bundle $\bee_p$ is 
a $GL(n,\C)$ bundle but not a $SL(n,\C)$-bundle, because it has 
non-trivial determinant. The obstruction class to lift the $PGL(n,\C)$ bundle ${\mathbb P}\bee$ 
to an $SL(n,\C)$ bundle is a class 
$B\in H^2(\M_{Dol}(SL(n,\C),\Z_n))\subset H^2(\M_{Dol}(SL(n,\C)),U(1))$, which gives us our $B$-field on $\M_{Dol}(SL(n,\C))$. Moreover $B$  
has \cite{HT4}[Section 3] 
a natural equivariant extension $\hat{B}\in H_{\Gamma}^2(\M_{Dol}(SL(n,\C)),U(1))$, giving us our $B$-field on $\M_{Dol}(PGL(n,\C))$. This 
$B$-field will come handy for our mirror symmetry discussions below.

For non-Abelian Hodge theory on a curve, see \cite{hitchin}, which gives
a gauge theoretical approach, and yields the natural hyperk\"ahler metrics
on our spaces. \cite{goldman-xia} introduces to the $GL(1,\C)$ case in detail.  
On the geometry and cohomology 
of ${\mathcal M}^1_{Dol}(SL(2,\C))$ see \cite{hausel2}.
\end{parag}
\vskip.5cm
\begin{parag}{The mixed Hodge structure  on  
non-Abelian Hodge cohomologies}
The main subject of this survey paper is the mixed Hodge polynomial 
of the (sometimes stringy, sometimes with a $B$-field) 
cohomology of the spaces $\M_{Dol}(G)$, $\M_{DR}(G)$ and 
$\M_{B}(G)$, 
for our three choices for $G= 
GL(n,\C), PGL(n,\C)$ or $SL(n,\C)$. As a notational convenience we may omit
$G$ and simply write $\M_{B}$, $\M_{DR}$ and $\M_{Dol}$, when it is clear what
$G$ should be, or if $G$ could be any of our groups. 

First take $G=GL(1,\C)$. From (\ref{gl1}) 
we can easily calculate the mixed Hodge polynomials as follows: 
\bes \begin{array}{c} H(x,y,t;\M_{B}(GL(1,\C)))=(1+xyt)^{2g}\\   
H(x,y,t;\M_{Dol}(GL(1,\C)))= H(x,y,t;\M_{DR}(GL(1,\C)))= 
(1+xt)^g(1+yt)^g.\label{glone}\end{array}\ees
It is remarkable that $H(x,y,t;\M_{B}(GL(1,\C)))\neq H(x,y,t;\M_{DR}(GL(1,\C)))$ even though the spaces are analytically isomorphic. Moreover 
we can explicitly see that the mixed Hodge structure on 
$H^k(\M_{Dol}(GL(1,\C)),\C)$ and 
$H^k(\M_{DR}(GL(1,\C)),\C)$ are pure, while on $H^k(\M_{B}(GL(1,\C)),\C)$ 
it is not. 

From a K\"unneth argument we also see that:
\bes H(x,y,t;\M_{Dol}(GL(n,\C)))= H(x,y,t;\M_{Dol}(PGL(n,\C))) H(x,y,t;\M_{Dol}(GL(1,\C))), \ees 
and similarly for the other two spaces. 
Thus the calculation for $GL(n,\C)$ is equivalent with the calculation 
for $PGL(n,\C)$. 

Now we list what we know about the cohomologies $H^*(\M,\C)$.
The Poincar\'e polynomials $P(t;{\mathcal M}^1(SL(2,\C)))$ and    
$P(t;{\mathcal M}^1(PGL(2,\C)))$ were calculated in \cite{hitchin}, while
the polynomials  $P(t;{\mathcal M}^1(SL(3,\C)))$ and    
$P(t;{\mathcal M}^1(PGL(3,\C)))$ have been calculated 
in \cite{gothen}. In both
papers Morse theory for a natural $\C^\times$ action on $\M_{Dol}$ was used (acting by multiplication of the Higgs field). The idea
is to calculate the Poincar\'e polynomial of the various fixed point components of this action, and then sum them up with a 
certain shift. The largest of the fixed point components,  
when $\phi=0$, 
is an important and well-studied space itself so we 
define it here as \begin{multline} \N(SL(n,\C)):=\\ \left\{ {\mbox{the moduli space of stable vector bundles
of fixed determinant bundle of degree {d}}}\label{smallguy}
\right\}.\end{multline} The Poincar\'e
polynomial of this space was calculated in \cite{harder-narasimhan} by 
arithmetic and in \cite{atiyah-bott} by gauge theoretical methods with
explicit formulas in \cite{zagier}. Thus its contribution to $P(t;\M(SL(n,\C)))$ is easy to handle. However the other components of the fixed point set of the natural circle action is  
 more cumbersome to determine already 
when $n=4$, and consequently this Morse theory approach 
has not been completed 
for $n\geq 4$.

As a running  example here we calculate from \cite{hitchin} 
the Poincar\'e polynomial of ${\mathcal M}^1(PGL(2,\C))$ when $g=3$:
\begin{multline} P(t;{\mathcal M}^1(PGL(2,\C)))=
\\=3\,{t}^{12}+12\,{t}^{11}+18\,{t}^{10}+32\,{t}^{9}+18\,{t}^{8}+12\,{t}^
{7}+17\,{t}^{6}+6\,{t}^{5}+2\,{t}^{4}+6\,{t}^{3}+{t}^{2}+1
\label{p2t}
\end{multline}

The cohomology ring of ${\mathcal M}^1_{Dol}(PGL(2,\C))$ has been described 
explicitly by
generators \cite{HT1} and relations \cite{HT2}. A result which 
proved to be essential to produce our main Conjecture~\ref{mainconj}. Finally
Markman \cite{markman} showed that for $PGL(n,\C)$ the universal
cohomology classes do generate the cohomology ring.  

Considering the mixed Hodge structure on the cohomology of our 
spaces the following result first appeared in \cite{mehta} using
a construction of \cite{HT1}. Here we present a simple proof. 

\begin{theorem} The mixed Hodge structure on $H^k(\M_{Dol},\C)$ is pure of weight $k$. \label{pure} 
\end{theorem}

\begin{proof} Recall the compactification 
$\overline{\mathcal M}^d_{Dol}$ of
 $\M_{Dol}$ constructed in \cite{hausel1}. 
From that paper  it follows
   that  $\overline{\mathcal M}^d_{Dol}$ is a projective orbifold, so its mixed
  Hodge structure on $H^k(\overline{\mathcal M}^d_{Dol},\C)$ is pure of weight $k$. Now
  \cite{hausel1} also implies that the natural map
  $H^*(\overline{\mathcal M}^d_{Dol},\C)\to H^*(\M_{Dol},\C)$ is 
surjective. The functoriality of mixed Hodge structures \cite{deligne2} completes the proof. 
\end{proof}

One can similarly prove the same result for $\M_{DR}$. 

\begin{theorem} The mixed Hodge structure on $H^k(\M_{DR},\C)$ is pure of
  weight $k$. 
\end{theorem} 

\begin{proof} As explained in \cite{HT4}[Theorem 6.2] one can 
deform the complex structure of $\overline{\mathcal M}^d_{Dol}$ to 
the projective orbifold
$\overline{\mathcal M}^d_{DR}$, which is the 
compactification of $\M_{DR}$ given by
Simpson in \cite{simpson5}. Now this way we see that 
the natural map 
$H^*(\overline{\mathcal M}^d_{DR},\C)\to H^*(\M_{DR},\C)$ is 
a surjection, getting our
result as in the previous proof.
\end{proof}

In fact the argument in \cite{HT4}[Theorem 6.2] shows that  
\begin{theorem}[HT4]  The mixed 
Hodge structure on $H^*(\M_{Dol},\C)$ 
is isomorphic with the mixed Hodge structure on 
$H^*(\M_{DR},\C)$. \label{dr=dol}
\end{theorem}

The mixed Hodge structure on $\M_{B}$ however has not 
been 
studied in the literature. We
will start the study of it later in this paper, where we will see
that this Hodge structure will in fact be very much not pure. 
But for now we explain our reason to be interested in these 
mixed Hodge structures on the spaces $\M_{Dol}$, $\M_{DR}$ and
$\M_{B}$. 
The reason is mirror symmetry:
\end{parag}

\section{Mirror symmetry conjectures}

Our starting point in \cite{HT4} was the observation that the pair
$\M_{DR}(SL(n,\C))$ together with the $B$-field $B^e$ and 
${\mathcal M}^d_{DR}(PGL(n,\C))$ with the $B$-field $\hat{B}^d$ 
 satisfy the 
geometric picture for mirror symmetry conjectured by 
Strominger--Yau--Zaslow \cite{SYZ} and modified for $B$-fields 
by Hitchin in \cite{hitchin3}. This geometric picture 
requires the existence of a special Lagrangian fibration on both spaces, so that the fibres are dual. In fact in \cite{HT4} it is shown
that the so-called Hitchin map \cite{hitchin2} provides the required special Lagrangian
fibration on our spaces, with dual Abelian varieties as fibers. For
details on this see \cite{HT4}[Section 3].

Our focus in this survey is the topological implications of this 
mirror symmetry. The following conjecture is what we call 
the topological 
mirror test for our SYZ-mirror partners. 

\begin{conjecture}[\cite{HT4}]\label{mirrorderham}
For all $d,e \in \Z$, satisfying $(d,n)=(e,n)=1$, we have
$$\Est^{B^e}\Left(x,y;\M_\DR(SL(n,\C))\Right) 
= \Est^{\hat B^d}\Left(x,y;{\mathcal M}^e_\DR(PGL(n,\C))\Right).$$

\end{conjecture}
\begin{conjecturemark} 
Since $M^d_\DR(SL(n,\C))$ is smooth, the left-hand side actually equals the $E$-polynomial
$E\Left(x,y;M^d_\DR(SL(n,\C))\Right)$, which is independent of $e$. This 
motivates the following:
\begin{conjecture}[HT4] For any two $d_1$ and $d_2$ as long as $(d_1,n)=(d_2,n)=1$ we have:  \beq E\Left(x,y;{\mathcal M}^{d_1}_{Dol}(SL(n,\C))\Right)=E\Left(x,y;{\mathcal M}^{d_2}_{Dol}(SL(n,\C))\Right).\eeq \label{coprime}
\end{conjecture}
This, if so, is quite interesting as the Betti numbers of 
  $\N(SL(n,\C))$, the moduli space of stable vector bundles, with
fixed determinant of degree $d$, which can be considered
as the ``K\"ahler version'' of $\M_{Dol}(SL(n,\C))$, is known to depend on 
$d$. 
Already
when $n=5$, Zagier's explicit formula \cite{zagier} 
for $P(t;{\mathcal N}^1(SL(5,\C)))$ and 
$P(t;{\mathcal N}^3(SL(5,\C)))$ are different. 
We will see strong support for this 
Conjecture~\ref{coprime} later in Corollary~\ref{betticoprime}.  
\end{conjecturemark}

\begin{conjecturemark} Conjecture~\ref{mirrorderham} 
was proved for $n=2$ and $n=3$ in \cite{HT4}. 
The proof proceeds by first transforming the calculation 
to $\M_{Dol}$ via Theorem~\ref{dr=dol} and then uses 
the Morse theoretic method of \cite{hitchin} and \cite{gothen}. 
It is unclear, however, how this method can be extended 
for $n\geq 4$.  
\end{conjecturemark}

\begin{conjecturemark} An important ingredient of the proofs was 
a modification of a result of Narasimhan--Ramanan 
\cite{narasimhan-ramanan} to Higgs bundles.  It 
describes the fixed points of the action by the elements of 
$Jac[n]$ on $\M_{Dol}(SL(n,\C))$. The fixed point sets will be
some lower rank $m|n$ Higgs moduli 
spaces $\M_{Dol}(SL(m,\C);\tilde{C})$ for a certain covering 
$\tilde{C}$ of $C$. Their cohomology enters in the stringy 
contribution to the right hand side of 
Conjecture~\ref{mirrorderham} (recall (\ref{bfieldstringy})). 
\end{conjecturemark}

\begin{parag}{Number theory to the rescue} 
Although our mirror symmetry  Conjecture~\ref{mirrorderham} 
is still open for $n\geq 4$, recently
some support for the validity of it has been achieved in 
form of progress 
on another related conjecture.

\begin{conjecture}[HRV] \label{mirrorbetti}
For all $d,e \in \Z$, such that $(d,n)=(e,n)=1$, we have
$$\Est^{B^e}\Left(x,y,\M_\B(SL(n,\C))\Right) 
= \Est^{\hat B^d}\Left(x,y,{\mathcal M}^e_\B(PGL(n,\C))\Right).$$
\end{conjecture}

This conjecture has been proved \cite{HRV} when $n$ is a prime 
and when $n=4$; which implies 
Conjecture~\ref{mirrorderham} on the level of Euler 
characteristic in these cases. 
The method of proof is arithmetic. We count 
rational
points of the variety ${\M_B}$ over a finite field $\F_q$, when
$n$ divides $q-1$, where $q=p^r$ is a prime power. Because this 
count will turn out to be a polynomial in $q$ we have the following
 result, which is basically the Weil conjecture for our special 
smooth affine varieties:

\begin{theorem}[\cite{HRV}] The $E$-polynomial of $\M_B$ has only $x^ky^k$ type terms, and the polynomial  $E(q):=E(\sqrt{q},\sqrt{q})$ agrees with the number of rational points of $\M_{B}(G)$ 
over $\F_q$:  \bes E(q)=\#(\M_{B}(G)(\F_q))  \ees
\end{theorem}

The count is then
possible because we only
need to count the solutions of the equation:
\bes [A_1,B_1]
\cdot \dots \cdot [A_g,B_g]=\xi_n, \ees
in the finite group of Lie type $G(\F_q)$, i.e.
so that $A_i,B_i\in G(\F_q)$, where $\xi_n\in G$ is a central element of order $n$. 

A simple modification of a theorem of Mednykh \cite{mednykh}, 
(which goes back to Frobenius--Schur \cite{frobenius-schur}, 
and has since  been obtained by many authors, 
for example by Freed--Quinn 
\cite{freed-quinn}[(5.19)]) then implies:

\begin{theorem} Let $G=SL(n,\C)$ or $G=GL(n,\C)$. Then the 
number of rational points on $\M_{B}(G)$ over a finite field 
$\F_q$, where $q=p^r$ is a prime power, 
with $n|(q-1)$  is given by the character  formula: 
\bes \#(\M_{B}(G)(\F_q))=\sum_{\chi \in {\rm Irr(G(\F_q))}} 
\frac{|G|^{2g-2}}{\chi(1)^{2g-1}} \chi(\xi_n) ,  \ees
where the sum is over all irreducible characters of the 
finite group $G(\F_q)$ of Lie type.
\end{theorem}

The two theorems above imply the following

\begin{corollary}[\cite{HRV}] The $E$-polynomial of the character variety
$\M_B(G)$ is given by the character formula:
\beq E(q)=\sum_{\chi \in {\rm Irr(G(\F_q))}} 
\frac{|G|^{2g-2}}{\chi(1)^{2g-1}} \chi(\xi_n).\eeq \label{eqcount}
\end{corollary} 

\begin{remark}  An immediate consequence of this formula is
the Betti analogue of Conjecture~\ref{coprime}. 
This follows from Corollary~\ref{eqcount} as that character
formula transforms by a Galois automorphism when one changes
from $d_1$ to $d_2$. Moreover because our varieties ${\mathcal M}^{d_1}_{B}(G)$ and  ${\mathcal M}^{d_2}_{B}(G)$  are Galois 
conjugate themselves, we can deduce 
that their Betti numbers agree, and 
presumably their mixed Hodge 
structure should also agree. In summary we have

\begin{corollary}[\cite{HRV}]  For all $d_1,d_2\in \Z$ \label{betticoprime}
as long as $(d_1,n)=(d_2,n)=1$ we have \beq E\Left(x,y;{\mathcal M}^{d_1}_B(G)\Right)=E\Left(x,y;{\mathcal M}^{d_2}_B(G)\Right)\eeq and 
\beq P\Left(t;{\mathcal M}^{d_1}_B(G)\Right)=
P\Left(t;{\mathcal M}^{d_2}_B(G)\Right).\eeq
\end{corollary}

Thus we get an affirmative answer for Conjecture~\ref{coprime} on the level of the Poincar\'e polynomials. 
In general Galois conjugate varieties tend to be (although need not 
be see e.g. \cite{serre})  homeomorphic over $\C$.   

\begin{problem} Are ${\mathcal M}^{d_1}_B(G\Right)$ and ${\mathcal M}^{d_2}_B(G\Right)$ homeomorphic for $(n,d_1)=(n,d_2)=1$?
\end{problem}
\end{remark}
\begin{remark}
In order to calculate the character formula in Corollary~\ref{eqcount}, we will need to know the 
values of irreducible characters of $G$ on central elements. 
Fortunately for $GL(n,\F_q)$ this has been calculated by Green \cite{green} and for $SL(n,\F_q)$ the required information, i.e. the value of the characters on central elements, was obtained
by Lehrer in \cite{lehrer}. In the next section we will show
an explicit result for the character formula for $GL(n,\F_q)$. 
\end{remark}

\begin{remark} Our mirror symmetry Conjecture~\ref{mirrorbetti} then can be 
translated to a complicated formula which is valid for the 
character tables of $PGL(n,\F_q)$ and $SL(n,\F_q)$. 
In particular we believe that by introducing punctures for our 
Riemann surfaces
a similar mirror symmetry conjecture would in fact capture the 
exact difference between the full character tables 
of $PGL(n,\F_q)$ 
and $SL(n,\F_q)$ (not just on central elements as above). 
This way our mirror symmetry proposal could be 
phrased as follows: {\em the differences between the 
character tables
of  $PGL(n,\F_q)$ and its Langlands dual $SL(n,\F_q)$ are 
governed by mirror symmetry}. In particular it is particularly
enjoyable
to see the effect of mirror symmetry on the differences between
the  character tables of 
$GL(2,\F_q)$ and $SL(2,\F_q)$, which were first 
calculated in 1907 by  Jordan 
\cite{jordan} and by Schur \cite{schur}.

\end{remark}
\end{parag}
\section{Explicit formulas for the E-polynomials} Here we 
calculate the $E$-polynomials of $\M_B(PGL(n,\C))$, which we 
denote by $E_n(q)$.  
We need to start with partitions.

 We write a
partition of $n$ as $\lambda =(\lambda_{1}\geq \lambda_{2}\geq \cdots \geq \lambda
_{l}>0)$, so that $\sum \lambda_i=n$.  The
{\it Ferrers diagram} $d(\lambda)$ 
of $\lambda $ is the set of lattice points
\begin{equation}\label{e:d(lambda)}\{(i,j)\in {\mathbb Z_{\leq 0}} \times
  {\mathbb N} : j <\lambda_{-i+1} \}.
\end{equation}
 The {\it arm length} $a(z)$ and
{\it leg length} $l(z)$ of a point $z\in d(\lambda )$ denote the number of points
strictly to the right of $z$ and below $z$, respectively, as indicated
in this example: 

 \begin{eqnarray*}\begin{array}[c]{cccccc}
 \bullet &      \bullet &       \bullet &
\bullet &       \bullet  \\

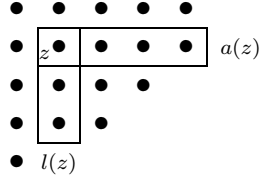
\cline{2-5}
\bullet &       \multicolumn{1}{|c|}{\llap{${}_{z}$}\bullet } & \bullet &
    \bullet&     \multicolumn{1}{c|}{\bullet }& {\scriptstyle a(z)} \\
\cline{2-5}
\bullet &       \multicolumn{1}{|c|}{\bullet }& \bullet &       \bullet \\

\bullet &       \multicolumn{1}{|c|}{\bullet }& \bullet \\
\cline{2-2}
\bullet &       \hbox to 0pt{\hss $\scriptstyle l(z)$\hss }\\
\end{array}
\end{eqnarray*}
where $\lambda=(5,5,4,3,1)$, $z=(-1,1)$, $a(z)=3$ and $l(z)=2$. The {\em hook length} then is defined as $$h(z)=l(z)+a(z)+1.$$

Let $$V_n(q)=E_n(q)  {q^{(1-g)n(n-1)}(q-1)^{2g-2}},$$ and
$$Z_n(q,T)=\exp\left(\sum_{r\geq 1} V_n (q^r) \frac{T^r}{r}\right).$$  We define
the Hook polynomials for a  partition $\lambda$ as follows
: 
$${\cal H}^{\lambda}(q)= 
\prod_{z \in d(\lambda)}
{q^{-l(z)}}(1-q^{h(z)}).$$

We can now formulate 
\begin{theorem}[\cite{HRV}] The $E$-polynomials of the character varieties
$\M_{B}(PGL(n,\C))$ for $n=1,2,3,\dots$ are given by the following 
generating function :
\begin{eqnarray}  
\prod_{n=1}^\infty Z_n(q,T^n) = \sum_{\lambda\in {\cal P}} 
({\cal H}^\lambda (q))^{2g-2} T^{|\lambda|},\end{eqnarray} 
where ${\cal P}$ is the set of all partitions. \label{enq} 
\end{theorem}

One simple corollary of this gives a new topological result:

\begin{corollary}[\cite{HRV}] The Euler characteristic of $\M(PGL(n,\C))$ is $\mu(n) n^{2g-3}$, where $\mu$ is the
  M\"obius function, i.e. $\mu(n)$ is the sum of primitive $n$th root of unities. \label{euler}
\end{corollary}

Another interesting application of the theorem is the following:

\begin{corollary} The $E$-polynomial $E_n(q)=E(q; \M_{B}(PGL(n,\C)))$
  is palindromic, i.e. it satisfies, what we call, the {\em curious}
  Poincar\'e duality:
$$q^{2N} E_n(1/q)=E_n(q),$$ where $2N=(n^2-1)(2g-2)$ is the complex
  dimension of $\M_{B}(PGL(n,\C))$. 
\end{corollary}

\begin{remark} In fact this result originates in the so-called
  Alvis--Curtis duality \cite{alvis,curtis} in the character theory of $GL(n,\F_{q})$,
  which is a duality between irreducible representations of
  $GL(n,\F_{q})$. In particular, if $\chi, \chi^\prime\in {\rm
  Irr}(GL(n,\F_q))$ are dual representations then the dimension 
  $\chi(1)$ is a polynomial in $q$ which satisfies
  $$q^{\frac{n(n-1)}{2}} \chi(1)(1/q)=\chi^\prime(1)(q).$$
\end{remark} 
For example when $n=2$ Theorem~\ref{enq} gives:
\beq E_2(q)=(q^2-1)^{2g-2}+q^{2g-2}(q^2-1)^{2g-2}-\frac{1}{2} q^{2g-2} (q-1)^{2g-2}-\frac{1}{2}q^{2g-2}(q+1)^{2g-2}\label{e2q},\eeq
when $g=3$ this gives 
\beq E(x,y;{\mathcal M}^1_{B}(PGL(2,\C)))= 
{q}^{12}-4\,{q}^{10}+6\,{q}^{8}-14\,{q}^{6}+6\,{q}^{4}-4\,{q}^{2}+1, 
\label{e2qg3}\eeq 
which is a palindromic polynomial indeed. Note also that there does not seem to be much in common with 
the Poincar\'e polynomial (\ref{p2t}). 

\section{A conjectured formula for mixed Hodge polynomials} 

Here we present the conjecture of \cite{HRV} on the $H$-polynomials of
the spaces $\M_B(PGL(n,\C))$.  
As usual 
we fix the curve $C$ and its genus $g$ and the group $PGL(n,\C)$ and
write $\M_B$ for $\M_B(PGL(n,\C))$ and $H_n(x,y,t)$ for $H(x,y,t;\M_B)$.   

Now we introduce rational functions $H_n(q,t)$ in two variables,
via a generating function. We let $$V_n(q,t)=H_n(q,t)  \frac{(qt^2)^{(1-g)n(n-1)}(qt+1)^{2g}}{(qt^2-1)(q-1)},$$ and
$$Z_n(q,t,T)=\exp\left(\sum_{r\geq 1} V_n (q^r,-(-t)^r) \frac{T^r}{r}\right).$$  We define
the $t$-deformed Hook polynomials for genus $g$ and partition $\lambda$ as follows: 
$${\cal H}_g^{\lambda}(q,t)= 
\prod_{x\in d(\lambda)}
\frac{(qt^2)^{(2-2g)l(x)}(1+q^{h(x)} t^{2l(x)+1})^{2g}}{(1-q^{h(x)} t^{2l(x)+2})(1-q^{h(x)}t^{2l(x)})}.$$
The following generating function then defines our rational functions
$H_n(q,t)$: \begin{eqnarray} \label{main} 
\prod_{n=1}^\infty Z_n(q,t,T^n) = \sum_{\lambda\in {\cal P}} {\cal H}^\lambda_g (q,t) T^{|\lambda|}.\end{eqnarray}
 Because for the character variety
we have that $h^{i,j;k}(\M_{B})=0$ provided that $i\neq j$ the following
conjecture describes $H_n(x,y,t)$ completely. 

\begin{conjecture}[\cite{HRV}] \label{mainconj}
 The mixed Hodge polynomials of the varieties
$\M_{B}(PGL(n,\C))$ are given by the generating function (\ref{main}): 
\bes H_n(\sqrt{q},\sqrt{q},t)=H_n(q,t)\ees
\end{conjecture}

Thus $H_n(q,t)$, which is a priori only a rational function, 
is conjectured to be the $H$-polynomial of the character
variety, so in the next conjecture we formalize  
our expectations from $H_n(q,t)$, with the addition of a curious, 
Poincar\'e duality-type of symmetry, which was 
in fact our most important guide to come up with these formulas:

\begin{conjecture} The rational functions $H_n(q,t)$ defined in the generating function of (\ref{main})
satisfy the following properties: 
\label{comb}
\begin{itemize}
\item $H_n(q,t)$ is a polynomial in $q$ and $t$.
\item Both the $q$ degree and the $t$ degree of the polynomial $H_n(q,t)$ agree with
$2N=2(n^2-1)(g-1)$. In fact the largest degree monomial in both variables is $(qt)^{2(n^2-1)(g-1)}$. 
\item All coefficients of $H_n(q,t)$ are non-negative integers. \label{combconj}
\item The coefficients of $H_n(q,t)=\sum h^i_j q^i t^j$ satisfy, what we call the {\em curious Poincar\'e duality}: 
\begin{eqnarray}h^{i-j}_{{N}-j}=h^{i+j}_{{N}+j} \label{curious} \end{eqnarray}
\end{itemize}
\end{conjecture}
In the following 
we list some checks and implications  of the above conjectures:
\vskip.4cm
\begin{conjecturemark} 
Computer calculations with Maple gives 
$H_n(q,t)$ from the above
generating function when $n=2,3,4$. In all these cases for small $g$ 
we do get a polynomial in $q$ and $t$ with
the expected degree and positive coefficients, satisfying the curious symmetry (\ref{curious}).
\end{conjecturemark}
\begin{conjecturemark} 
Using the explicit description of the cohomology ring of
$\M_{B}$ in 
\cite{HT2}, one can write
down a monomial basis for $H^*(\M_{B},\C)$ in the tautological
generators. 
Then in turn one can
figure out the action of the Frobenius on these generators, which in turn provides a 
formula for the mixed Hodge polynomial 
of $\M_{B}$. This formula can be brought to the form 
\begin{eqnarray} H_2(\sqrt{q},\sqrt{q},t)&=&
\frac{(q^2t^3+1)^{2g}}{(q^2t^2-1)(q^2t^4-1)}+
\frac{q^{2g-2}t^{4g-4}(q^2t+1)^{2g}}{(q^2-1)(q^2t^2-1)}-\nonumber \\
&-&\frac{1}{2}\frac{q^{2g-2}t^{4g-4}(qt+1)^{2g}}{(qt^2-1)(q-1)}-\frac{1}{2}
\frac{ q^{2g-2}t^{4g-4}(qt-1)^{2g}}{(q+1)(qt^2+1)},
\label{h2qt}\end{eqnarray} 
which agrees with the 
conjectured one through (\ref{main}), and clearly reduces to (\ref{enq}), 
when $t=-1$. For example when $g=3$, this gives
\begin{multline}H(\sqrt{q},\sqrt{q},t;\M(PGL(2,\C)))=\\ 
= {t}^{12}{q}^{12}+{t}^{12}{q}^{10}+6\,{t}^{11}{q}^{10}+{t}^{12}{q}^{8}+
{t}^{10}{q}^{10}+6\,{t}^{11}{q}^{8}+16\,{t}^{10}{q}^{8}+6\,{t}^{9}{q}^
{8}+{t}^{10}{q}^{6}+{t}^{8}{q}^{8}+26\,{t}^{9}{q}^{6}+\\ +16\,{t}^{8}{q}^{
6}+6\,{t}^{7}{q}^{6}+{t}^{8}{q}^{4}+{t}^{6}{q}^{6}+6\,{t}^{7}{q}^{4}+
16\,{t}^{6}{q}^{4}+6\,{t}^{5}{q}^{4}+{t}^{4}{q}^{4}+{t}^{4}{q}^{2}+6\,
{t}^{3}{q}^{2}+{t}^{2}{q}^{2}+1, \label{h2qtg3}
\end{multline}
which is a common refinement of (\ref{p2t}) when $q=1$ and of 
$(\ref{e2qg3})$
when $t=-1$. Note also how the curious Poincar\'e duality appears
when one refines the Poincar\'e polynomial (\ref{p2t}), which does not possess
any kind of symmetry, to the mixed Hodge
polynomial (\ref{h2qtg3}).

\end{conjecturemark}
\begin{conjecturemark} Note that $P_{n}(t)=H_{n}(1,t)$ should be the Poincar\'e polynomial
of the character 
variety, which is the same as the Poincar\'e
polynomial of the diffeomorphic Higgs moduli space $\M_{Dol}$. 
For $n=2$ Hitchin in \cite{hitchin} calculated the Poincar\'e
polynomial 
of this latter space, and an easy calculation shows that if one substitutes $q=1$ into
(\ref{h2qt}) we get $P_2(t)=H_2(1,t)$, the Poincar\'e polynomial 
of Hitchin. For $n=3$ Gothen in \cite{gothen} calculated $P_3(t)$. 
Because it is so pleasant to work with a formula like (\ref{h2qt}), we also give 
what our Conjecture~\ref{main} gives in the $n=3$ case:  
\begin{eqnarray*} H_3(q,t)&=&{\frac { \left( {q}^{3}{t}^{5}+1 \right)
      ^{2\,g} 
\left( {q}^{2}{t}^{3}+1 \right) ^{2\,g}}{ \left( {q}^{3}{t}^{6}-1 \right)  \left( {q}^{3}{t
}^{4}-1 \right)  \left( {q}^{2}{t}^{4}-1 \right)  \left( {q}^{2}{t}^{2
}-1 \right) }}+{\frac {{q}^{6\,g-6}{t}^{12\,g-12} \left( {q}^{3}t+1
 \right) ^{2\,g} \left( {q}^{2}t+1 \right) ^{2\,g}}{ \left( {q}^{3}{t}
^{2}-1 \right)  \left( {q}^{3}-1 \right)  \left( {q}^{2}{t}^{2}-1
 \right)  \left( {q}^{2}-1 \right) }}+\\
&+&{\frac {{q}^{4\,g-4}{t}^{8\,g-8}
 \left( {q}^{3}{t}^{3}+1 \right) ^{2\,g} \left( qt+1 \right) ^{2\,g}}{
 \left( {q}^{3}{t}^{4}-1 \right)  \left( {q}^{3}{t}^{2}-1 \right) 
 \left( q{t}^{2}-1 \right)  \left( q-1 \right) }}+\frac{1}{3}\,{\frac {{q}^{6
\,g-6}{t}^{12\,g-12} \left(  \left( qt+1 \right) ^{2\,g} \right) ^{2}}
{ \left( q{t}^{2}-1 \right) ^{2} \left( q-1 \right) ^{2}}}-\\ &-&\frac{1}{3}\,{
\frac {{q}^{6\,g-6}{t}^{12\,g-12} \left( {q}^{2}{t}^{2}-qt+1 \right) ^
{2\,g}}{ \left( {q}^{2}{t}^{4}+q{t}^{2}+1 \right)  \left( {q}^{2}+q+1
 \right) }}-{\frac {{q}^{4\,g-4}{t}^{8\,g-8} \left( {q}^{2}{t}^{3}+1
 \right) ^{2\,g} \left( qt+1 \right) ^{2\,g}}{ \left( {q}^{2}{t}^{4}-1
 \right)  \left( {q}^{2}{t}^{2}-1 \right)  \left( q{t}^{2}-1 \right) 
 \left( q-1 \right) }}-\\&-&{\frac {{q}^{6\,g-6}{t}^{12\,g-12} \left( {q}^{
2}t+1 \right) ^{2\,g} \left( qt+1 \right) ^{2\,g}}{ \left( {q}^{2}{t}^
{2}-1 \right)  \left( {q}^{2}-1 \right)  \left( q{t}^{2}-1 \right) 
 \left( q-1 \right) }}.\end{eqnarray*}  
Now it is a nice exercise to show that $H_3(1,t)$ does indeed produce 
(the corrected version\footnote{One accidental
mistake in the calculation of 
\cite{gothen} were pointed out in (10.3) of  
\cite{HT4}.} of) Gothen's complicated looking formula in \cite{gothen}. 

It is also worth noting that many of the individual terms in
$H_n(q,t)$ have poles at 
$q=1$, however
according to our conjecture these poles somehow cancel each other. 
\end{conjecturemark}

\begin{conjecturemark}
 When $g=0$, we know from the definitions that $H_1(x,y,t)=1$ and
$H_n(x,y,t)=0$ otherwise. 
One can deduce the same
from Conjecture~\ref{mainconj} by applying
part f of Theorem 2.10 in \cite{garsia-haiman} to calculate the right
hand side of (\ref{main}). Moreover Conjecture~\ref{comb} has the 
same flavour as the main conjecture in \cite{garsia-haiman} about
$q,t$ Catalan numbers, which was in turn proved
by Haiman in \cite{haiman} 
using some subtle intersection theory on the 
Hilbert scheme of $n$ points on $\C^2$. Apart from the fact that this
Hilbert scheme is also a hyperk\"ahler manifold, the similarities 
between the two conjectures are rather surprising. 
\end{conjecturemark}

\begin{conjecturemark}
\noindent When $g=1$ we have $H_n(x,y,t)=1$ for every $n$, but  
this we
could not prove from (\ref{main}) for $H_n(q,t)$. 
\end{conjecturemark}

\begin{conjecturemark} 
Let us look at the conjecture (\ref{curious}). Recall that $H^2$ 
of our varieties are exactly one dimensional, generated by a class, call it 
$[\omega]$, which is the K\"ahler class in the complex structure of
$\M_{Dol}$.
 This carries the weight $q^2t^2$ in the 
mixed Hodge structure. We have the following hard Lefschetz
type of conjecture which enhances the curious 
Poincar\'e duality of the conjecture of (\ref{curious}):

\begin{conjecture} If $L$ denotes the map by multiplication with
  $[\omega]$,
 then we conjecture that the map 
$$ L^k: H^{{N}-k,N-k;i-k}(\M_{B}(PGL(n,\C)))\to H^{N+k,N+k;i+k}(\M_{B}(PGL(n,\C)))$$ is an isomorphism.  
\end{conjecture} 

Interestingly this conjecture implies a theorem of \cite{hausel3} 
that the Lefschetz map $L^k:H^{N-k}\to H^{N+k}$ is injective for
$\M_{Dol}$, 
and it is explained 
there how this weak version of Hard Lefschetz, when applied to toric
hyperk\"ahler varieties, yields
new inequalities for the $h$-numbers of matroids. See also
\cite{hausel-sturmfels} for the original argument on toric hyperk\"ahler 
varieties. 
Furthermore this conjecture 
can be proved when $n=2$ using the explicit description of the cohomology ring
in \cite{HT2}. The general case  can also be thought of as an analogue
of the Faber 
conjecture \cite{faber} on the cohomology of the moduli space of curves,
which is another non-compact variety whose cohomology ring is conjectured 
to satisfy a certain form of the Hard Lefschetz theorem. 
\end{conjecturemark}

\begin{conjecturemark}
There are two subspaces of the cohomology $H^*(\M_{B},\C)$
 which are
 particularly interesting. 
One is the middle dimensional cohomology $H^{2N}(\M_{B},\C)$,
 which is the
 top non-vanishing cohomology. 
The mixed Hodge structure will break it into parts with respect to the $q$-degree. The
curious Poincar\'e dual (\ref{curious}) of these spaces are also
 interesting: it is
 easy to see that they are exactly the pure part of the mixed Hodge structure i.e. spaces of the form $H^{i,i;2i}$.  (Another significance of the pure part is that if there is a smooth projective
compactification of the variety  then its image is in this pure part.) 
Thus it would already be interesting  to get the pure part of $H_n(q,t)$. In fact it is easy to identify the pure part in our 
case with what we call the Pure ring, which is the subring of $H^*(\M,\C)$ generated by the  
tautological classes $a_i\in H^{2i}(\M,\C)$ for $i=2,..,n$ 
(the other
 tautological classes, which generate the cohomology ring, 
are not pure classes).  

For example when $n=2$, it was known 
\cite{hitchin}
that the middle degree cohomology of the Higgs moduli space
$\M_{Dol}(PGL(2,\C))$ 
is $g$ dimensional. 
The Pure ring was determined in \cite{HT2}, 
and it was found to be
$g$ dimensional due to the relation $\beta^g=0$ (where 
$\beta=a_2$). 
Thus these two seemingly
unrelated observations are dual to each other via our curious Poincar\'e duality (\ref{curious}). To see this curious duality in action let us
recall the formula (\ref{h2qtg3}). The terms which contain the top degree
$12$ in $t$ are $t^{12}q^{12}$, $t^{12} q^{10}$ and $t^{12} q^8$, which
are curious Poincar\'e dual via (\ref{curious}) to the terms $1$, $t^4q^2$
and $t^8q^4$, which is exactly the ring generated by the degree four class
$\beta$, which has additive basis $1$, $\beta$ and $\beta^2$.

The analogous 
ring, generated by the corresponding classes 
$a_2,\dots,a_n\in H^*(\N,\C)$, 
which a priori is a quotient of our Pure ring 
(as $\N\subset \M_{Dol}$ naturally), was studied for
the moduli space $\N$ of rank $n$, degree $d$ stable bundles (with
 $(n,d)=1$) in \cite{earl-kirwan}, 
where they in particular found the top non-vanishing degree of this 
ring to be $2n(n-1)(g-1)$. 
Computer calculations for our conjecture for $n=2,3,4$ also show 
that 
our conjectured Pure ring has the same $1$-dimensional top 
degree. This and the known situation for $n=2$ (see \cite{HT2}), 
yields the following 

\begin{conjecture} \label{pontrequ} The Pure rings of 
$\M_{Dol}$
 and $\N$, i.e. the subrings of the cohomology rings generated
by the classes $a_2,\dots,a_n$  
are isomorphic. In particular, unlike the whole cohomology ring of
 $\N$, it does not depend on $d$.  
\end{conjecture}

Now we explain a 
combinatorial consequence of this conjecture. First
we extract  a conjectured formula for $PP_n(t)$ 
the Poincar\'e polynomial  of the Pure ring. Indeed we only have
to deal with monomials in Conjecture~\ref{mainconj} whose 
$t$-degree
is double of their $q$-degree. 

Thus let $$PV_n(t)=PP_n(t)\frac{ t^{2(1-g)n(n-1)}}{(t^2-1)},$$ and
$$PZ_n(t,T)=\exp\left(\sum_{r\geq 1} PV_n (t^r) \frac{T^r}{r}\right).$$ We now  define
the pure part of the $t$-deformed Hook polynomials for genus 
$g$ and
partition $\lambda$ as follows:
\begin{eqnarray*}{\cal P H}^{\lambda}_{g}(t)=
t^{4(1-g)n(\lambda^\prime)} \prod_{x\in d(\lambda); a(x)=0}
\frac{1}{(1-t^{2 h(x)})},\end{eqnarray*}
where $$n(\lambda^\prime):=\sum_{z\in d(\lambda)} l(z).$$
Thus we get the conjecture that $PP_n(t)$ is given by
\begin{eqnarray} \label{pontr}
\prod_{n=1}^\infty PZ_n(t,T^n) = \sum_{\lambda\in {\cal P}} {\cal P
  H}^\lambda_g (t) T^{|\lambda|}.\end{eqnarray}

Now combining the two conjectures above we can formulate a new conjecture:

\begin{conjecture} The rational functions $PP_{n}(t)$ defined in (\ref{pontr}) satisfy 
  \begin{itemize} \item $PP_{n}(t)$ is a polynomial in $t$
 \item all coefficients of $PP_{n}(t)$ are non-negative integers
 \item The degree of $PP_{n}(t)$ is $2n(n-1)(g-1)$, and the coefficients of the leading term is $1$
\end{itemize}
 
\end{conjecture}
    
So for example when $n=3$ our conjecture gives for the Poincar\'e polynomial of the Pure ring:
$$PP_3(t)={\frac {1}{ \left( {t}^{6}-1 \right)  \left( {t}^{4}-1 \right) }}+t^{12\,g-12} -{
\frac {{t}^{8\,g-8}}{t^2-1}}+\frac{1}{3}\,{\frac {{t}^{12\,g-12}}{ \left( t^2-1
 \right) ^{2}}}-\frac{1}{3}\,{\frac {{t}^{12\,g-12}}{{t}^{4}+t^2+1}}-{\frac {{t}^{8\,g-8}}{ \left( {t}^{4}-1
 \right)  \left( t^2-1 \right) }}+{\frac {{t}^{
12\,g-12}}{t^2-1}}$$
\end{conjecturemark}

\begin{conjecturemark} Interestingly we can modify the formula of Conjecture~\ref{mainconj} to get a
   conjectured formula for the mixed Hodge polynomial of
   $\M_{Dol}$. Recall from Theorem~\ref{pure} that the
   mixed Hodge structure on $H^k(\M_{Dol},\C)$ is pure of
   weight $k$, thus this mixed Hodge polynomial is equivalent 
with the $E$-polynomial.

We now introduce polynomials $H_n(q,x,y)$ of three variables. 
Let $$V_n(q,x,y)=H_n(q,x,y)  \frac{(qxy)^{(1-g)n(n-1)}(qx+1)^{g} (qy+1)^{g}}{(qxy-1)(q-1)},$$ and
$$Z_n(q,x,y,T)=\exp\left(\sum_{r\geq 1} V_n (q^r,-(-x)^r,-(-y)^r) \frac{T^r}{r}\right).$$  We define
the $(x,y)$-deformed Hook polynomials for genus $g$ and partition $\lambda$ as follows: 
\begin{eqnarray}{\cal H}_g^{\lambda}(q,x,y)= 
\prod_{z\in d(\lambda)}
\frac{(qxy)^{(2-2g)l(z)}(1+q^{h(z)} y^{l(z)}
  x^{l(z)+1})^{g}(1+q^{h(z)} x^{l(z)} y^{l(z)+1})^{g}} {(1-q^{h(z)}
  (xy)^{l(z)+1})(1-q^{h(z)}(xy)^{l(z)})}.\label{xydefhook}\end{eqnarray} The following generating
function defines $H_n(q,x,y)$:
\begin{eqnarray}  
\prod_{n=1}^\infty Z_n(q,x,y,T^n) = \sum_{\lambda\in {\cal P}} {\cal H}^\lambda_g (q,x,y) T^{|\lambda|}.\end{eqnarray}

Clearly we have $H_n(q,t,t)=H_n(q,t)$ which says that a specialization
of $H_n(q,x,y)$ gives the mixed Hodge polynomial $H_n(q,t)$ of
$\M_B$.  The following conjecture says that another
specialization gives the mixed Hodge polynomial of $\M_{Dol}$ and
$\M_{DR}$.

\begin{conjecture} $H_n(q,x,y)$ is a polynomial with non-negative
 integer coefficients
with specialization 
$H_n(1,x,y)$ equal to the $E$-polynomial
  of  the Higgs moduli space  $\M_{Dol}(PGL(n,\C))$. \label{higgs}
\end{conjecture} 

Thus we have a mysterious formula $H_n(q,x,y)$ which
 specializes, on one hand to the $H$-polynomial of the character
 variety, and on the other hand 
to the mixed Hodge polynomial of the Higgs (or
 equivalently flat connection) moduli space. It would be very
 interesting to find a geometrical meaning for $H_n(q,x,y)$. 

Checks on this Conjecture~\ref{higgs} include a proof for $n=2$ and
$n=3$, (one can easily modify Hitchin's and Gothen's argument to get 
the Hodge polynomial instead of the Poincar\'e polynomial of the Higgs
moduli space) and
also computer checks that the shape of the polynomial $H_n(1,x,y)$ is the
expected one when $n=4$. 

Consider now the specification $H_n(q,-1,y)$. Interestingly, the 
corresponding specification of the $(x,y)$-deformed Hook 
polynomials (\ref{xydefhook}) becomes a polynomial, showing that 
at least $H_n(q,-1,y)$ is 
a polynomial. We get an even nicer formula if we make the further 
specification $H_n(1,-1,y)$ which by Conjecture~\ref{higgs} 
should be the Hirzebruch $y$-genus
of the moduli space of Higgs bundles $\M_{Dol}$. Namely, for $g>1$, most of the $(x,y)$
deformed Hook polynomials vanish, when one substitutes first $x=-1$ and
then $q=1$. Indeed, the only partitions which will have a non-zero 
contribution to the $y$-genus are the partitions of the form
 $n=1+1+\dots+1$; when $l(z)=0$ only for once. 
This in turn
gives the following closed formula for the conjectured $y$-genus of 
$\M_{Dol}$:

\begin{conjecture} The Hirzebruch $y$-genus of $\M_{Dol}(PGL(n,\C))$, for $g>1$, 
equals $$ (1-y+\dots +(-y)^{n-1})^{g-1}\sum_{m | n} \frac{\mu(m)}{m} \left( (-y)^{n(n-n/m)}m
\prod_{i=1}^{n/m-1} (1-(-y)^{mi})^2 \right)^{g-1}$$ 
\end{conjecture}

In particular note that the term corresponding to $m=1$, is exactly
the known $y$-genus of $\N$ (see \cite{narasimhan-ramanan}). The rest
thus should be thought as contribution of the other fixed point 
components of the circle action on $\M_{Dol}$. Of course this conjectured
formula gives the known specialization of Corollary~\ref{euler} 
at $y=-1$, while the $y=1$ specialization gives $\mu(n) n^{g-2}$ when
$n$ is odd, and $0$ when $n$ is even. The specialization at $y=1$  
can be 
thought of as the signature of the pairing on the rationalized 
circle equivariant 
cohomology of $\M_{Dol}$ as defined in \cite{hausel-proudfoot}. 
\end{conjecturemark}

\begin{conjecturemark} Finally we discuss how to obtain a conjecture 
for the mixed Hodge polynomial of $\M_{B}(SL(n,\C))$. For the mixed Hodge polynomial of  
$\M_{Dol}(SL(n,\C))$ the mirror symmetry conjecture~\ref{mirrorderham}, 
together with Conjecture~\ref{higgs} imply a conjecture. For $\M_{B}(SL(n,\C))$ the mixed Hodge polynomial contains more information than the $E$-polynomial. In order thus to have a conjecture on 
$H_n(x,y,t;\M_{B}(SL(n,\C)))$ 
a mirror symmetry conjecture is needed on the level of the $H$-polynomial.
We finish by formulating such a conjecture, generalizing 
Conjecture~\ref{mirrorbetti} for $H$-polynomials:

\begin{conjecture} For all $d,e \in \Z$, with $(d,n)=(e,n)=1$ we have
$$\Hst^{B^e}\Left(x,y,t;\M_\B(SL(n,\C))\Right) 
= \Hst^{\hat B^d}\Left(x,y,t;{\mathcal M}^e_\B(PGL(n,\C))\Right),$$
where $\Hst^B$ is the stringy mixed Hodge polynomial twisted with a 
$B$-field, which can be defined identically as $\Est^B$ is
defined in (\ref{bfieldstringy}).
\end{conjecture}

\end{conjecturemark}

\end{document}